\documentclass[10pt,a4paper]{amsart}

\usepackage{setspace}
\usepackage[utf8x]{inputenc}
\usepackage[T1]{fontenc}

\usepackage{helvet}

\usepackage{amsthm}
\usepackage{amsmath}
\usepackage{verbatim}
\usepackage{amssymb} 
\usepackage{resizegather}
\usepackage{tikz}
\usepackage{tikz-cd}
\usepackage{babel} 
\usepackage{mathrsfs}
\usepackage{graphicx}
\usepackage{epstopdf}
\usepackage{epigraph}

\usepackage{pgfplots}
\usepgfplotslibrary{fillbetween}
\pgfplotsset{compat=1.13}

\usepackage[bottom]{footmisc}

\usepackage{hyperref}

\newtheorem{theorem}{Theorem}[section]

\newtheorem{corollary}{Corollary}[theorem]

\newtheorem{definition}{Definition}

\title{On Isolated Real Singularities II}
\author{Lars Andersen}

\begin{document}

\maketitle
\begin{abstract} In this article we apply the results in \cite{Lars} to real $ADE$-singularities. In effect, we show in Corollary \hyperref[ADE]{\ref*{ADE}} how \cite[Theorem 4.1]{Lars} and \cite[Theorem 5.1]{Lars}, taken in tandem, enables us to find in the case of an $ADE$-singularity, the Poincaré polynomials of the Milnor fibres, except for the case with a $D_k^{-, s}$-singularity (see the list \hyperref[list]{\ref*{list}} below) where $k$ is even.
\end{abstract}
\begin{center}
    \text{Classification: }\textbf{14-XX}
\end{center}

\section{Introduction}
In this paper we shall consider isolated singularities $f: \mathbb{R}^{n+1}\to\mathbb{R}$ of the following form:
\begin{enumerate}\label{list}
\item $A_k^{\pm, s}: f=x^{k+1}\pm y^2+\sum_{i=1}^{t} x_i^2-\sum_{j=t+1}^{s+t} x_j^2,\quad k\geq 2$.
\item $D_k^{\pm, s}: f=x^2 y\pm y^{k-1}+\sum_{i=1}^{t} x_i^2-\sum_{j=t+1}^{s+t} x_j^2,\quad k\geq 4$
\item $E_6^{\pm s}: f=x^3 \pm y^4+\sum_{i=1}^{t} x_i^2-\sum_{j=t+1}^{s+t} x_j^2$
\item $E_7^s: f=x^3+xy^3+\sum_{i=1}^{t} x_i^2-\sum_{j=t+1}^{s+t} x_j^2$
\item $E_8^s: f=x^3+y^5+\sum_{i=1}^{t} x_i^2-\sum_{j=t+1}^{s+t} x_j^2$
\end{enumerate}
where $t+s=n-1, s,t\geq 0$ and where $(x,y, x_1,\dots, x_{n-1})$ are coordinates on $\mathbb{R}^{n+1}$. To simplify the notation we shall write $\mathcal{F}^{+}=\mathcal{F}_{\eta}^{+}(f)$ for a positive Milnor fibre and $\mathcal{F}^{-}=\mathcal{F}^{-}(f)$ for a negative Milnor fibre at the origin. 

\subsection*{Acknowledgments} The author wish to express his sincere gratitude to the\\
Laboratoire de Mathématiques at USMB and especially to his thesis supervisors Georges Comte and Michel Raibaut. The material formed in this article and its successor formed a part of the authors thesis. On behalf of all authors, the corresponding author states that there is no conflict of interest.
\subsection*{Data Sharing} Data sharing not applicable to this article as no datasets were generated or analysed during the current study.

\section{On the Topology of Real Curve $ADE$-Singularities}
The idea will be to first treat the case of curves (namely $n=1$) and find for each of the classes above a $\mathbb{R}$-morsification of a particularly simple character, with few critical points. It is then easy to find $\mathbb{R}$-morsifications in the general case (namely $n\geq 2$).\\

The following result gives the Poincaré polynomials of the Milnor fibres of $ADE$-singularities of curves.

\begin{theorem}\label{ADE thing} For a topological space $X$ let $\beta(X)$ denote the Poincaré polynomial in singular homology with $\mathbb{Z}$-coefficients. Then 
\begin{enumerate}
\item $A_k^{+}: f(x,y)=x^{k+1}+y^2$: 
$$\beta(\bar{\mathcal{F}}^{+})=\beta(\bar{\mathcal{F}}^{-})=1,\qquad\text{if }k\equiv 0\quad(\text{mod 2})$$
$$\beta(\bar{\mathcal{F}}^{+})=1+u,\quad \beta(\bar{\mathcal{F}}^{-})=0\qquad\text{if }k\equiv 1\quad(\text{mod 2}).$$

 \item $A_k^{-}: f(x,y)=x^{k+1}-y^2$: \[\beta(\bar{\mathcal{F}}^{+})=\beta(\bar{\mathcal{F}}^{-})= \left\{
    \begin{array}{ll}
      1,\qquad k\equiv 0\quad(\text{mod 2})\\
      2,\qquad  k\equiv 1 \quad(\text{mod 2})
\end{array} \right. \]
\item $D_k^{+}: f(x,y)=x^2y+y^{k-1},\quad k\geq 4$, \[\beta(\bar{\mathcal{F}}^{+})=\beta(\bar{\mathcal{F}}^{-})= \left\{
    \begin{array}{ll}
      1,\qquad k\equiv 0\quad(\text{mod 2})\\
      2,\qquad  k\equiv 1 \quad(\text{mod 2})
\end{array} \right. \]
\item $D_k^{-}: f(x,y)=x^2y-y^{k-1},\quad k\geq 4$, \[\beta(\bar{\mathcal{F}}^{\pm})= \left\{
    \begin{array}{ll}
      3,\qquad k\equiv 0\quad(\text{mod 2})\\
      2,\qquad  k\equiv 1 \quad(\text{mod 2})
\end{array} \right. \]

\item $E_6^{\pm}: f(x,y)=x^3\pm y^4$, $\beta(\bar{\mathcal{F}}^{+})=\beta(\bar{\mathcal{F}}^{-})=1.$
\item $E_7: f(x,y)=x^3+xy^3$, $\beta(\bar{\mathcal{F}}^{+})=\beta(\bar{\mathcal{F}}^{-})=2.$
\item $E_8: f(x,y)=x^3+y^5,$ $\beta(\bar{\mathcal{F}}^{+})=\beta(\bar{\mathcal{F}}^{-})=1.$
\end{enumerate} 
\end{theorem}

\begin{proof} We are going to find, in each of the cases above, a $\mathbb{R}$-morsification having either one unique critical point, or no critical points at all and in the case $D_k^{-}$ with $k$ odd, a $\mathbb{R}$-morsification with two critical points. To begin with we treat the case of an $A_k$-singularity.

\begin{enumerate}
\item Consider $f: \mathbb{R}^2\to\mathbb{R}$ given by $f(x,y)=x^{k+1}+y^2$.\\
Suppose that $k\equiv 0\quad(\text{mod }2)$ and let 
$$F:\mathbb{R}^2\times[0, 1]\to\mathbb{R}$$
$$F(x,y,t)=x^{k+1}+tx+y^2.$$
Then 
$$\text{Jac}(f_t)=[(k+1)x^{k}+t, 2y],$$
so $f_t$ has no critical point for $t\neq 0$. Hence applying \cite[Theorem 1]{Lars} we get 
$$\beta(\bar{\mathcal{F}}^{+})=\beta(\bar{\mathcal{F}}^{-})=1$$
since the Milnor  fibres are contractible.\\
Suppose that $k\equiv 1\quad(\text{mod }2)$. Then $\bar{\mathcal{F}}^{-}=\emptyset$ because $f$ is nonnegative. Letting 
$$F:\mathbb{R}^2\times[0, 1]\to\mathbb{R},\qquad F(x,y,t)=x^{k+1}-(k+1)tx+y^2$$
then 
$$\text{Jac}(f_t)=[(k+1)(x^{k}-t), 2y],$$
so that there is a unique critical point $p=(t^{1/k}, 0)$ and then 
\[ \text{Hess}(f_t)(p)=
\begin{pmatrix}
    k(k+1) t^{(k-1)/k}  & 0 \\
    0  & 2
\end{pmatrix}
.\]
Since $t>0$ it follows that $\lambda(p)=0$. Therefore we can apply \cite[Corollary 5.1.2]{Lars} to deduce that 
$$\beta(\bar{\mathcal{F}}^{+})=1+u.$$
\item Consider $f: \mathbb{R}^2\to\mathbb{R}$ given by $f(x,y)=x^{k+1}-y^2$.\\
Suppose that $k\equiv 0\quad(\text{mod }2)$ and define 
$$F:\mathbb{R}^2\times[0, 1]\to\mathbb{R},\qquad F(x,y,t)=x^{k+1}+tx-y^2.$$
Then 
$$\text{Jac}(f_t)=[(k+1)x^{k}+t, -2y],\qquad t\in [0,1]$$
Since $k$ is even, $f_t$ has no critical points whenever $t$ is nonzero. Therefore the positive and negative Milnor fibres are both contractible, by \cite[Theorem 4.1]{Lars} so 
$$\beta(\bar{\mathcal{F}}^{\pm})=1.$$
Suppose that $k\equiv 1\quad(\text{mod }2)$. Then one considers $F:\mathbb{R}^2\times[-1, 0]\to\mathbb{R}$ given by the same formula as in the case $k$ even and one gets instead that $f_t=F(x,t)$ has a unique critical point 
$$p=((-t/(k+1))^{1/k}, 0).$$ 
in which the Hessian matrix is
 \[ \text{Hess}(f_t)(p)=
\begin{pmatrix}
    \frac{k(k+1)}{(k+1)^{(k-1)/k}} (-t)^{(k-1)/k}  & 0 \\
    0  & -2
\end{pmatrix}
.\]
Since $t$ is negative there is only one negative eigenvalue so the index is $\lambda(p)=1$. One applies \cite[Corollary 5.1.2]{Lars} and concludes that 
$$\beta(\bar{\mathcal{F}}^{+})=\beta(\bar{\mathcal{F}}^{-})=2.$$

\item Let $k\geq 4$ and consider a $D_k^{+}$-singularity $f(x,y)=x^2y+y^{k-1}$.\\ 
Suppose first that $k\equiv 0\quad(\text{mod }2)$ and consider 
$$F:\mathbb{R}^2\times[-1, 0]\to\mathbb{R},$$
$$F(x,y,t)=x^2y+2tx^2+y^{k-1}-ty.$$
Then 
$$\text{Jac}(f_t)=[2x(y+2t), x^2+(k-1)y^{k-2}-t]$$
has nonmaximal rank in a point $p=(x, y)$ if and only if either one of the following sets of equations is satisfied
 \[\left\{
    \begin{array}{ll}
      x=0,\quad (k-1)y^{k-2}=t,\quad\text{or}\\
      x^2=t\left(1-2^{k-2}(k-1)t^{k-3}\right),\quad y=-2t.
\end{array} \right. \]
Remark first of all that if $t\neq 0$ then $(k-1)y^{k-2}=t$ has no real solutions because $k$ is even and $t<0$. Furthermore, since $t<0$ and since $(k-3)$ is odd,
$$\forall t<0,\qquad t^{k-3}<0<\frac{1}{(k-1)2^{k-2}}$$
hence $1-2^{k-2}(k-1)t^{k-3}>0$. Therefore the second set of equations has no real solutions for $t\neq 0$ either. As a consequence $f_t: \mathbb{R}^2\to \mathbb{R}$ has no critical points whenever $t\neq 0$. Thus applying \cite[Theorem 4.1]{Lars} one gets that the Milnor fibres are contractible and so $\beta(\bar{\mathcal{F}}^{+})=\beta(\bar{\mathcal{F}}^{+})=1$ by homotopy invariance.\\
Suppose now that $k\equiv 1\quad(\text{mod }2)$. Put 
$$I=(-\frac{1}{\left((k-1) 2^{k-2}\right)^{1/(k-3)}}, 0]$$
and consider the $\mathbb{R}$-morsification
$$F:\mathbb{R}^2\times I\to\mathbb{R},$$
$$F(x,y,t)=x^2y+2tx^2+y^{k-1}-ty.$$
Since $k-3\equiv 0\quad(\text{mod }2)$ one has 
$$\forall t\in I,\qquad t^{k-3}<\frac{1}{(k-1) 2^{k-2}}$$
Therefore if $t\in I\setminus\{0\}$ then 
$$x^2=t\left(1-2^{k-2}(k-1)t^{k-3}\right)$$
has no real solutions. As a consequence the only critical point of $f_t$ is $p(t)=(0, (t/k-1)^{k-2})$. The Hessian $\text{Hess}(f_t)(p(t))$ is 
\[
\begin{pmatrix}
    2t(2+\frac{t^{k-3}}{(k-1)^{k-2}})  & 0 \\
    0  & (k-1)(k-2)(\frac{t}{k-1})^{(k-2)(k-3)}
\end{pmatrix}
.\]
Since 
$$k-3\equiv (k-2)(k-3)\equiv 0\quad (\text{mod }2)$$
it follows that if $t\in I\setminus\{0\}$ then $\lambda(p)=1$. Using \cite[Corollary 5.1.2]{Lars} one deduces that $$\beta(\bar{\mathcal{F}}^{+})=\beta(\bar{\mathcal{F}}^{-})=2.$$

\item Let $k\geq 4$ and consider a $D_k^{-}$-singularity $f(x,y)=x^2y-y^{k-1}$.\\
Suppose that $k\equiv 1\quad(\text{mod }2)$ and let 
$$F:\mathbb{R}^2\times[0, 1]\to\mathbb{R},\qquad F(x,y,t)=x^2y-y^{k-1}+(k-1)ty.$$
Then 
$$\text{Jac}(f_t)=[2xy, x^2+(k-1)(t-y^{k-2})].$$
If $t\neq 0$ then there is only one critical point 
 $$p(t)=(0, t^{1/(k-2)})$$
and one find the Hessian to be given by
\[ \text{Hess}(f_t)(p(t))=
\begin{pmatrix}
    2t^{1/(k-2)}  & 0 \\
    0  & -(k-1)(k-2)t^{(k-3)/(k-2)}
\end{pmatrix}
.\] 
Therefore $\lambda(p(t))=1$ whenever $t\neq 0$ hence
$$\beta(\bar{\mathcal{F}}^{+})=\beta(\bar{\mathcal{F}}^{-})=\beta(\mathbb{S}^0)=2$$
by \cite[Corollary 5.1.2]{Lars}.\\
Suppose that $k\equiv 0\quad(\text{mod }2)$ and let 
$$F:\mathbb{R}^2\times[0, 1]\to\mathbb{R},\qquad F(x,y,t)=x^2y-y^{k-1}+(k-1)ty.$$
Then 
$$\text{Jac}(f_t)=[2xy, x^2+(k-1)(t-y^{k-2})].$$
If $t\neq 0$ then there are two critical points 
$$p_1(t)=(0, t^{1/(k-2)}),\qquad p_2(t)=(0, -t^{1/(k-2)})$$
and the Hessian matrices are 
\[ \text{Hess}(f_t)(p(t))=
\begin{pmatrix}
    \pm 2t^{1/(k-2)}  & 0 \\
    0  & \mp (k-1)(k-2)t^{(k-3)/(k-2)}
\end{pmatrix}
.\] 

As a consequence, since $t>0$, the Morse indices are 
$$\lambda(p_1)=\lambda(p_2)=1.$$
Since $n+1-\lambda(p_i)=\lambda(p_i)$ one can apply \cite[Theorem 4.1]{Lars} to deduce that 
$$\beta(\bar{\mathcal{F}}_{\eta}^{\pm})=\beta(\mathbb{S}^0\vee\mathbb{S}^0)=3.$$

\item Consider a $E_6^{\pm}$-singularity $f^{\pm}(x,y)=x^3\pm y^4$. If 
$$F^{\pm}:\mathbb{R}^2\times[0,1]\to\mathbb{R},\qquad F^{\pm}(x,y,t)=x^3+ 3tx\pm y^4$$
and if $f_t^{\pm}=F^{\pm}(\cdot, t)$ then 
$$\text{Jac}(f_t^+)=[3x^2+3t, 4y^3],$$
$$\text{Jac}(f_t^{-})=[3x^2+3t, -4y^3]$$
so that neither $f_t^{+}$ nor $f_t^{-}$ has any critical points whenever $t\neq 0$ hence 
$$\beta(\bar{\mathcal{F}}^{+})=\beta(\bar{\mathcal{F}}^{-})=1$$
by \cite[Theorem 4.1]{Lars}. The same holds for $E_8$-singularities because if one defines 
$$F:\mathbb{R}^2\times[0,1]\to\mathbb{R},\qquad F(x,y,t)=F(x,y,t)=x^3+3tx+y^5$$
then 
$$\text{Jac}(f_t)=[3x^2+3t, 5y^4]$$
so $f_t$ has no critical points whenever $t$ is nonzero hence the Milnor fibres are contractible.
\item It remains the case of an $E_7$-singularity. In this case, put  
$$F:\mathbb{R}^2\times [0,1]\to\mathbb{R},\quad F(x,y,t)=x^3+3tx+xy^3+ty^3.$$
The Jacobian matrix is 
$$\text{Jac}(f_t)=[3(x^2+t)+y^3, 3y^2(x+t)]$$
which has nonmaximal rank in $p=(x,y)$ if and only if either 
\[\left\{
    \begin{array}{ll}
      x^2+t=0,\quad y=0,\quad\text{or}\\
      x=-t,\quad y^3=-3(t^2+t).
\end{array} \right. \]
However if $t\neq 0$ then the first set of equations has no real solutions and one finds that $f_t$ has a unique critical point in
$$p(t)=(-t, -(3)^{1/3}(t^2+t)^{1/3})$$
The Hessian matrix is given by
\[ \begin{pmatrix}
6x & 3y^2\\
3y^2 & 6y(x+t)
\end{pmatrix}
\]
and one finds that evaluated in $p(t)$ its eigenvalues are 
$$\lambda_{+}(t)=-3t+ \sqrt{9t^2+3^{10/3}(t^2+t)^{4/3}},$$
$$\lambda_{+}(t)=-3t-\sqrt{9t^2+3^{10/3}(t^2+t)^{4/3}}.$$
For $t$ sufficiently small, $\lambda_{+}(p(t))>0$ and $\lambda_{-}(p(t))<0$. Hence the Morse index at $p=p(t)$ is $\lambda(p)=1$ and, applying \cite[Corollary 5.1.2]{Lars} one gets that 
$$\beta(\bar{\mathcal{F}}^{+})=\beta(\bar{\mathcal{F}}^{-})=2.$$
\end{enumerate}

\end{proof}

\section{On the Topology of Real$ADE$-Singularities of Higher Codimension}
From this we can easily deduce the corresponding result for higher codimensions, except for the case $D_k^{-}$ with $k$ even and $n\geq 2$.

\begin{corollary}\label{ADE} The Poincaré polynomials in singular homology with integer coefficients of the Milnor fibres of the isolated singularities $A_k^{\pm s}, D_k^{+ s}$ and $E_i^s$ are given as follows.
\begin{enumerate}
 \item $A_k^{-s}, D_k^{+ s}$: \[\beta(\bar{\mathcal{F}}^{+})= \left\{
    \begin{array}{ll}
      1,\qquad k\equiv 0\quad(\text{mod 2})\\
      1+u^{n-s-1},\qquad  k\equiv 1 \quad(\text{mod 2})
\end{array} \right.\]
\[\beta(\bar{\mathcal{F}}^{-})= \left\{
    \begin{array}{ll}
      1,\qquad k\equiv 0\quad(\text{mod 2})\\
      1+u^{s},\qquad  k\equiv 1 \quad(\text{mod 2})
\end{array} \right.\]
\item $E_7^s:$ 
$$\beta(\bar{\mathcal{F}}^{+})=1+u^{n-s-1},$$
$$\beta(\bar{\mathcal{F}}^{-})=1+u^{s}.$$

\item $E_6^{\pm s}, E_8^s:$ $\beta(\bar{\mathcal{F}}^{+})=\beta(\bar{\mathcal{F}}^{-})=1.$
\item $A_k^{+ s}:$
\[\left\{
    \begin{array}{ll}
      \beta(\bar{\mathcal{F}}^{+})=1+u^{n-s},\qquad k\equiv 1\quad(\text{mod 2})\\
      \beta(\bar{\mathcal{F}}^{-})=1+u^{s-1},\qquad  s\neq 0, k\equiv 1 \quad(\text{mod 2})\\
      \beta(\bar{\mathcal{F}}^{-})=0,\qquad  s=0, k\equiv 1 \quad(\text{mod 2})\\
      \beta(\bar{\mathcal{F}}^{+})=\beta(\bar{\mathcal{F}}^{-})=1,\qquad k\equiv 0\quad(\text{mod }2).
\end{array} \right. \]
\end{enumerate}
The Poincaré polynomials of the Milnor fibres of $D_k^{- s}$-singularities for $k$ odd and $n>1$ are 
\[D_k^{- s}: \left\{
    \begin{array}{ll}
      \beta(\bar{\mathcal{F}}^{+})=1+u^{n-s-1},\qquad k\equiv 1\quad(\text{mod 2})\\
      \beta(\bar{\mathcal{F}}^{-})=1+u^{s},\qquad  k\equiv 1 \quad(\text{mod 2})
\end{array} \right.\]
 
\end{corollary}

\begin{proof} Consider an $\mathbb{R}$-morsification 
$$\tilde{F}: \mathbb{R}^2\times\bar{U}\to\mathbb{R},\qquad f_t=\tilde{F}(\cdot,t)$$
of any of the $ADE$-singularities $\tilde{f}: \mathbb{R}^2\to\mathbb{R}$ given in the proof of Theorem \hyperref[ADE thing]{\ref*{ADE thing}} and let  
$$f:\mathbb{R}^{n+1}\to\mathbb{R},\qquad f=\tilde{f}(x,y)+\sum_{i=1}^t x_i^2-\sum_{j=t+1}^{s+t} x_j^2.$$
If 
$$F: \mathbb{R}^{n+1}\times \bar{U}\to\mathbb{R},\qquad F=\tilde{F}+\sum_{i=1}^t x_i^2-\sum_{j=t+1}^{s+t} x_j^2$$
then the critical points of $f_t=F(\cdot, t)$ are the points $p=(\tilde{p},0,\dots,0)$ with $\tilde{p}$ a critical point of $\tilde{f}_t$. Moreover these are Morse whenever $\tilde{p}$ is Morse, with Hessian 

\[
\text{Hess}(f_t)(p)= \begin{pmatrix}
\text{Hess}(\tilde{f}_t)(\tilde{p}) & 0 & 0\\
0 & 2I_t & 0 \\
0 & 0 & -2I_s
\end{pmatrix}
\]

where $I_t$ and $I_s$ denotes the identity matrices of size $t\times t$ and $s\times s$, respectively. As a consequence, if $\lambda(\tilde{p})$ is the Morse index of $\tilde{f}_t$ at $\tilde{p}$ then $\lambda(p)=s+\lambda(\tilde{p})$ is the Morse index of $f_t$ at $p$. Since in each of the cases there is either one critical point or none at all the result follows immediately from \cite[Corollary 5.1.2]{Lars} and from \cite[Theorem 4.1]{Lars}.

\end{proof}

\section{Concluding Remarks}
We can remark here that G. Fichou (see \cite{fichou2008}) has classified the real $ADE$-singularities above using the blow-Nash equivalence relation. The following definitions can be found on page 184 of \cite{Arn}. 
\begin{definition}[{\cite[II]{Arn}}] Let $X$ be a manifold and suppose that a Lie group $G$ acts on $X$. Let $x\in X$. The \emph{modality} of $x$ under the action of $G$ is the least number $m$ such that there exists a neighborhood $x\in U\subset M$ and a covering of $U$ by finitely many $m$-parameter families of orbits. If $m=0$ then one says that $x$ is \emph{simple}. If $m=1$ one says that $x$ is \emph{unimodal} and if $m=2$ one says that $x$ is \emph{bimodal}.
\end{definition}

In other words, a point $x\in X$ is simple if it has a neighborhood $U$ intersecting only finitely many orbits of $G$.  

\begin{definition}[{\cite[II]{Arn}}] Let $f: (\mathbb{R}^{n+1},0)\to (\mathbb{R},0)$ be a real analytic function germ. The \emph{modality} of $f$ is the modality of its jet $j^k f$, for $k\in\mathbb{N}$ sufficiently large, in the space of jets of functions $\mathbb{R}^{n+1}\to \mathbb{R}$ having a critical point $0\in \mathbb{R}^{n+1}$ and a critical value $0\in\mathbb{R}$, under the action of the Lie group $\text{Diff}(\mathbb{R}^{n+1})$ of diffeomorphisms of $\mathbb{R}^{n+1}$. 
\end{definition}

Returning to Fichou's result, he showed that if 
$$f, g: (\mathbb{R}^{n+1}, 0)\to (\mathbb{R},0)$$ 
are two germs of Nash functions with one of them simple (as a germ of real analytic function), say $f$, then $f, g$ are blow-Nash equivalent if and only if they are analytically equivalent, and then $g$ is simple as well. So for simple Nash function germs the blow-Nash equivalence coincides with the analytic equivalence. Consequently, the list \hyperref[list]{\ref*{list}} given in the beginning of this section gives a complete list of simple singularities up to the blow-Nash equivalence.

\bibliographystyle{plain}
\bibliography{biblio.bib}

\end{document}